\documentclass[11pt,twoside]{article}

\usepackage[latin1]{inputenc}
\usepackage{epsfig}
\usepackage{amsfonts}
\usepackage{cite}
\usepackage{color}
\definecolor{lgrey}{gray}{0.99}

\title{
  {\huge Complex Analysis of Real Functions \\[1.5ex]}
  IV: Non-Integrable Real Functions }

\author{
  \Large Jorge L. deLyra\footnote{Email: delyra@latt.if.usp.br} \\
  Department of Mathematical Physics \\
  Physics Institute \\
  University of São Paulo }

\date{May 28, 2017}

\newcommand{\ii}{\mbox{\boldmath$\imath$}}

\newcommand{\e}[1]{\,{\rm e}^{#1}}

\newcommand{\ldot}{\mbox{\Large$\cdot$}\!}

\newtheorem{definition}{Definition}

\newtheorem{theorem}{Theorem}

\newtheorem{proof}{Proof}[theorem]
\newcommand{\Colon}{{\hspace{-0.4em}\bf:}}

\begin{document}\maketitle

\vspace{-1.0ex}
\begin{abstract}
  \noindent
  In the context of the complex-analytic structure within the unit disk
  centered at the origin of the complex plane, that was presented in a
  previous paper, we show that a certain class of non-integrable real
  functions can be represented within that same structure. In previous
  papers it was shown that essentially all integrable real functions, as
  well as all singular Schwartz distributions, can be represented within
  that same complex-analytic structure. The large class of non-integrable
  real functions which we analyze here can therefore be represented side
  by side with those other real objects, thus allowing all these objects
  to be treated in a unified way.
\end{abstract}

\section{Introduction}\label{Sec01}

In a previous paper~\cite{CAoRFI} we introduced a certain complex-analytic
structure within the unit disk of the complex plane, and showed that one
can represent essentially all integrable real functions within that
structure. In another, subsequent previous paper~\cite{CAoRFII}, we showed
that one can also represent within the same structure the singular objects
known as distributions, loosely in the sense of the Schwartz theory of
distributions~\cite{DTSchwartz}, which are also known as generalized real
functions. In this paper we will show that a large class of non-integrable
real functions can also be represented within that same structure. All
these objects will be interpreted as parts of this larger complex-analytic
structure, within which they can be treated and manipulated in a robust
and unified way.

The construction presented in~\cite{CAoRFI}, which leads to the inclusion
of all integrable real functions in the complex-analytic structure, starts
from the Fourier coefficients of an arbitrarily given real function on the
unit circle. Since these coefficients are defined as integrals involving
the real function, there is the obvious necessity that these real
functions be integrable on that circle. Therefore, that construction will
not work for functions which fail to be integrable there. However, once we
are in the space of inner analytic functions within the open unit disk,
the concept of integral-differential chains of inner analytic functions,
which was introduced in~\cite{CAoRFI}, comes to our rescue. It will allow
us to extend the representation within the complex-analytic structure to a
large class of non-integrable real functions.

For ease of reference, we include here a one-page synopsis of the
complex-analytic structure introduced in~\cite{CAoRFI}. It consists of
certain elements within complex analysis~\cite{CVchurchill}, as well as of
their main properties.

\paragraph{Synopsis:} The Complex-Analytic Structure\\

\noindent
An {\em inner analytic function} $w(z)$ is simply a complex function which
is analytic within the open unit disk. An inner analytic function that has
the additional property that $w(0)=0$ is a {\em proper inner analytic
  function}. The {\em angular derivative} of an inner analytic function is
defined by

\noindent
\begin{equation}
  w^{\ldot}(z)
  =
  \ii
  z\,
  \frac{dw(z)}{dz}.
\end{equation}

\noindent
By construction we have that $w^{\ldot}(0)=0$, for all $w(z)$. The {\em
  angular primitive} of an inner analytic function is defined by

\begin{equation}
  w^{-1\ldot}(z)
  =
  -\ii
  \int_{0}^{z}dz'\,
  \frac{w(z')-w(0)}{z'}.
\end{equation}

\noindent
By construction we have that $w^{-1\ldot}(0)=0$, for all $w(z)$. In terms
of a system of polar coordinates $(\rho,\theta)$ on the complex plane,
these two analytic operations are equivalent to differentiation and
integration with respect to $\theta$, taken at constant $\rho$. These two
operations stay within the space of inner analytic functions, they also
stay within the space of proper inner analytic functions, and they are the
inverses of one another. Using these operations, and starting from any
proper inner analytic function $w^{0\ldot}(z)$, one constructs an infinite
{\em integral-differential chain} of proper inner analytic functions,

\begin{equation}
  \left\{
    \ldots,
    w^{-3\ldot}(z),
    w^{-2\ldot}(z),
    w^{-1\ldot}(z),
    w^{0\ldot}(z),
    w^{1\ldot}(z),
    w^{2\ldot}(z),
    w^{3\ldot}(z),
    \ldots\;
  \right\}.
\end{equation}

\noindent
Two different such integral-differential chains cannot ever intersect each
other. There is a {\em single} integral-differential chain of proper inner
analytic functions which is a constant chain, namely the null chain, in
which all members are the null function $w(z)\equiv 0$.

A general scheme for the classification of all possible singularities of
inner analytic functions is established. A singularity of an inner
analytic function $w(z)$ at a point $z_{1}$ on the unit circle is a {\em
  soft singularity} if the limit of $w(z)$ to that point exists and is
finite. Otherwise, it is a {\em hard singularity}. Angular integration
takes soft singularities to other soft singularities, and angular
differentiation takes hard singularities to other hard singularities.

Gradations of softness and hardness are then established. A hard
singularity that becomes a soft one by means of a single angular
integration is a {\em borderline hard} singularity, with degree of
hardness zero. The {\em degree of softness} of a soft singularity is the
number of angular differentiations that result in a borderline hard
singularity, and the {\em degree of hardness} of a hard singularity is the
number of angular integrations that result in a borderline hard
singularity. Singularities which are either soft or borderline hard are
integrable ones. Hard singularities which are not borderline hard are
non-integrable ones.

Given an integrable real function $f(\theta)$ on the unit circle, one can
construct from it a unique corresponding inner analytic function $w(z)$.
Real functions are obtained through the $\rho\to 1_{(-)}$ limit of the
real and imaginary parts of each such inner analytic function and, in
particular, the real function $f(\theta)$ is obtained from the real part
of $w(z)$ in this limit. The pair of real functions obtained from the real
and imaginary parts of one and the same inner analytic function are said
to be mutually Fourier-conjugate real functions.

Singularities of real functions can be classified in a way which is
analogous to the corresponding complex classification. Integrable real
functions are typically associated with inner analytic functions that have
singularities which are either soft or at most borderline hard. This ends
our synopsis.

\vspace{2.6ex}

\noindent
When we discuss real functions in this paper, some properties will be
globally assumed for these functions, just as was done in the previous
papers~\cite{CAoRFI,CAoRFII,CAoRFIII} leading to this one. These are
rather weak conditions to be imposed on these functions, that will be in
force throughout this paper. It is to be understood, without any need for
further comment, that these conditions are valid whenever real functions
appear in the arguments. These weak conditions certainly hold for any real
functions that are obtained as restrictions of corresponding inner
analytic functions to the unit circle.

The most basic condition is that the real functions must be measurable in
the sense of Lebesgue, with the usual Lebesgue
measure~\cite{RARudin,RARoyden}. The second global condition we will
impose is that the functions have no removable singularities. The third
and last global condition is that the number of hard singularities on the
unit circle be finite, and hence that they be all isolated from one
another. There will be no limitation on the number of soft singularities.

The material contained in this paper is a development, reorganization and
extension of some of the material found, sometimes still in rather
rudimentary form, in the
papers~\cite{FTotCPI,FTotCPII,FTotCPIII,FTotCPIV,FTotCPV}.

\section{Preliminary Considerations}\label{Sec02}

Here we will quickly review some results of previous papers in this
series, which are necessary for the work currently at hand. We will also
establish some equally necessary groundwork for the discussion that
follows in Sections~\ref{Sec03} and~\ref{Sec04}.

\subsection{Singularities of Real Functions}\label{Ssec0201}

In this paper we will discuss non-integrable real functions, which
therefore have on the unit circle hard singularities which are not just
borderline hard ones. It is therefore necessary to consider in greater
detail the issue of the singularities of real functions. Therefore, let us
discuss the translation of the classification of singularities, which has
been established in~\cite{CAoRFI} for inner analytic functions in the unit
disk of the complex plane, to the case of real functions on the unit
circle.

First of all, let us discuss the concept of a removable singularity. This
is a well-known concept for analytic functions in the complex plane. What
we mean by a removable singularity in the case of real functions on the
unit circle is a singular point such that both lateral limits of the
function to that point exist and result in the same finite real value, but
where the function has been arbitrarily defined to have some other finite
real value. This is therefore a point were the function can be redefined
by continuity, resulting in a continuous function at that point. These
are, therefore, trivial singularities, which we will simply rule out of
our discussions in this paper.

The concept of a singularity is the same, namely a point where the
function is not analytic. The concepts of soft and hard singularities are
carried in a straightforward way from the case of complex functions to
that of real functions. The only significant difference is that the
concept of the limit of the function to a point is now taken to be the
real one, along the unit circle. The existence of the limit is a weaker
condition in the real case, because in the complex case the limit must
exist on and be independent of the continuum of directions along which one
can take it along the complex plane, while in the real case the limit only
has to exist and be the same along the two lateral directions. This means
that, if one takes a path across a singularity of an analytic function in
order to define a real function, the singularity of the real function may
be soft even if the singularity of the analytic function is hard, and the
singularity of the real function may be integrable even if the singularity
of the analytic function is not.

At this point it is interesting to note that it might be useful to
consider classifying all inner analytic functions $w(z)=u(\rho,\theta)+\ii
v(\rho,\theta)$ according to whether or not they have the same basic
analytic properties, taken in the complex sense, as the corresponding real
functions $u(1,\theta)$, taken in the real sense. One could say that a
{\em regular} inner analytic function is one that has, at all its singular
points on the unit circle, the same status as the corresponding real
function, regarding the most fundamental analytic properties. By contrast,
an {\em irregular} inner analytic function would be one that fails to have
the same status as the corresponding real function, regarding one or more
of the these analytic properties, such as those of integrability,
continuity, and a given level of multiple differentiability.

The gradations corresponding to soft and hard singularities can be
implemented for real functions in terms of the integrability or
non-integrability of the singularities. The soft singularities are all
integrable, and a singularity which is both hard and integrable is
necessarily a borderline hard singularity, that can be immediately
associated to the degree of hardness zero. The degree of softness of an
isolated soft singularity is the number of differentiations with respect
to $\theta$ of the real function $f(\theta)$, in a neighborhood of the
singular point, which are required for the singularity to become a
borderline hard one. Observing that a soft singularity of degree one means
that the function is continuous but not differentiable, an alternative
definition is that the degree of softness is one plus the number of
differentiations with respect to $\theta$ of the real function
$f(\theta)$, in a neighborhood of the singular point, which are required
for the function to become non-differentiable at that point. The remaining
problem is that of associating a degree of hardness to non-integrable hard
singularities. This is the most important part of this structure in regard
to our work in this paper, so let us examine it in more detail.

Let us assume that the real function $f(\theta)$ has an {\em isolated}
non-integrable singular point at $\theta_{1}$, which is therefore a hard
singularity. What we mean here by this singularity being isolated is that
the function has no other non-integrable hard singularities in a
neighborhood of that point, and is thus integrable in each one of the two
sides of that neighborhood. Consider then two points within that
neighborhood on the unit circle, say $\theta_{\ominus}$ and
$\theta_{\oplus}$, one to the left and the other to the right of
$\theta_{1}$, so that we have
$\theta_{\ominus}<\theta_{1}<\theta_{\oplus}$. It follows, therefore, that
the function $f(\theta)$ is an integrable real function on the two closed
intervals $[\theta_{\ominus},\theta_{1}-\varepsilon_{\ominus}]$ and
$[\theta_{1}+\varepsilon_{\oplus},\theta_{\oplus}]$, where
$\varepsilon_{\ominus}$ and $\varepsilon_{\oplus}$ are any sufficiently
small strictly positive real numbers, such that
$\theta_{\ominus}<\theta_{1}-\varepsilon_{\ominus}$ and
$\theta_{1}+\varepsilon_{\oplus}<\theta_{\oplus}$. Therefore, we can
integrate the function $f(\theta)$ within these two intervals, starting at
some arbitrary reference points $\theta_{0,\ominus}$ and
$\theta_{0,\oplus}$ strictly within each interval, defining in this way a
pair of sectional primitives $f_{\ominus}^{-1\prime}(\theta)$ and
$f_{\oplus}^{-1\prime}(\theta)$, one within each interval,

\noindent
\begin{eqnarray}
  f_{\ominus}^{-1\prime}(\theta)
  & = &
  \int_{\theta_{0,\ominus}}^{\theta}d\theta'\,
  f(\theta'),
  \nonumber\\
  f_{\oplus}^{-1\prime}(\theta)
  & = &
  \int_{\theta_{0,\oplus}}^{\theta}d\theta'\,
  f(\theta'),
\end{eqnarray}

\noindent
where we have that both the argument $\theta$ and the reference points
$\theta_{0,\ominus}$ and $\theta_{0,\oplus}$ are within the corresponding
closed integration intervals. Therefore, on the left interval we have that
$\theta_{\ominus}\leq\theta\leq\theta_{1}-\varepsilon_{\ominus}$ and that
$\theta_{\ominus}\leq\theta_{0,\ominus}\leq\theta_{1}-\varepsilon_{\ominus}$,
while on the right interval we have that
$\theta_{1}+\varepsilon_{\oplus}\leq\theta\leq\theta_{\oplus}$ and that
$\theta_{1}+\varepsilon_{\oplus}\leq\theta_{0,\oplus}\leq\theta_{\oplus}$. Note
that, since $f(\theta)$ is integrable on those closed intervals, it
follows that the piecewise primitive $f^{-1\prime}(\theta)$ formed by the
pair of sectional primitives $f_{\ominus}^{-1\prime}(\theta)$ and
$f_{\oplus}^{-1\prime}(\theta)$ is limited on those same closed intervals,
and therefore is also integrable on them. Therefore, this process of
sectional integration of the functions can be iterated indefinitely. Let
us assume that we iterate this process $n$ times, thus obtaining the
$n$-fold primitive $f^{-n\prime}(\theta)$.

In general one cannot take the limits $\varepsilon_{\ominus}\to 0$ or
$\varepsilon_{\oplus}\to 0$, because the singularity at $\theta_{1}$ is a
non-integrable one and therefore the limits of the asymptotic integrals on
either side will diverge. However, if there is a number $n$ of successive
sectional integrations such that the resulting primitive
$f^{-n\prime}(\theta)$ has a borderline hard singularity at $\theta_{1}$,
being therefore integrable on the whole closed interval
$[\theta_{\ominus},\theta_{\oplus}]$, then the two limits
$\varepsilon_{\ominus}\to 0$ and $\varepsilon_{\oplus}\to 0$ can both be
taken. In this case we may say that the hard singularity of $f(\theta)$ at
$\theta_{1}$ has degree of hardness $n$. In Section~\ref{Sec04} we will
have the opportunity to use these ideas regarding sectional integration on
closed intervals that avoid non-integrable singular points.

\subsection{Piecewise Polynomial Real Functions}\label{Ssec0202}

In~\cite{CAoRFII} we introduced the concept of piecewise polynomial real
functions, in the context of the discussion of singular Schwartz
distributions. These are the integrable real functions that are obtained
on the integration side of the integral-differential chains to which the
singular distributions belong, the latter being on the differentiation
side of the same chains, with respect to the central position of the delta
``function'' itself. This concept will return to the discussion in this
paper, where it will play an essential role. Therefore, let us review the
definitions and the relevant results established in that previous paper.

Let us assume that we have a real function which is defined by polynomials
in sections of the unit circle. These sections are open intervals
separated by a finite set of points which are singularities of that real
function. Let there be $N\geq 1$ singular points
$\{\theta_{1},\ldots,\theta_{N}\}$. It follows that in this case we must
separate the unit circle into a set of $N$ contiguous sections, consisting
of open intervals between the singular points, that can be represented as
the sequence

\begin{equation}\label{EQSeqSec}
  \left\{
    \rule{0em}{2ex}
    (\theta_{1},\theta_{2}),
    \ldots,
    (\theta_{i-1},\theta_{i}),
    (\theta_{i},\theta_{i+1}),
    \ldots,
    (\theta_{N},\theta_{1})
  \right\},
\end{equation}

\noindent
where we see that the sequence goes around the unit circle, and where we
adopt the convention that each section $(\theta_{i},\theta_{i+1})$ is
numbered after the singular point $\theta_{i}$ at its left end.

In~\cite{CAoRFII} we established a general notation for these piecewise
polynomial real functions, as well as a formal definition for them. Here
is the definition: given a real function $f_{(n)}(\theta)$ that is defined
in a piecewise fashion by polynomials in $N\geq 1$ sections of the unit
circle, with the exclusion of a finite non-empty set of $N$ singular
points $\theta_{i}$, with $i\in\{1,\ldots,N\}$, so that the polynomial
$P_{i}^{(n_{i})}(\theta)$ at the $i^{\rm th}$ section has order $n_{i}$,
we denote the function by

\begin{equation}\label{EQPieFun}
  f_{(n)}(\theta)
  =
  \left\{
    P_{i}^{(n_{i})}(\theta),
    i\in\{1,\ldots,N\}
  \right\},
\end{equation}

\noindent
where $n$ is the largest order among all the $N$ orders $n_{i}$. We say
that $f_{(n)}(\theta)$ is a {\em piecewise polynomial real function} of
order $n$. Note that, being made out of finite sections of polynomials,
the real function $f_{(n)}(\theta)$ is always an integrable real function.

These functions have some important properties, which were established
in~\cite{CAoRFII}. First and foremost, if the integrable real function
$f_{(n)}(\theta)$ is a non-zero piecewise polynomial function of order
$n$, then and only then its derivative $f_{(n)}^{(n+1)\prime}(\theta)$
with respect to $\theta$ is a superposition of a non-empty set of delta
``functions'' and derivatives of delta ``functions'' on the unit circle,
with their singularities located at some of the points $\theta_{i}$, and
of nothing else. In particular, the derivative
$f_{(n)}^{(n+1)\prime}(\theta)$ is identically zero within all the open
intervals defining the sections. However, this derivative cannot be
identically zero over the whole unit circle. In fact, it is impossible to
have a non-zero piecewise polynomial real function of order $n$, such as
the one described above, that is also continuous and differentiable to the
order $n$ on the whole unit circle. It follows that $f_{(n)}(\theta)$ can
be globally differentiable at most to order $n-1$, so that its $n^{\rm
  th}$ derivative is a discontinuous function, and therefore its
$(n+1)^{\rm th}$ derivative already gives rise to singular distributions.

In short, every real function that is piecewise polynomial on the unit
circle, of order $n$, when differentiated $n+1$ times, so that it becomes
zero within the open intervals corresponding to the existing sections,
will always result in the superposition of some non-empty set of singular
distributions with their singularities located at the points between two
consecutive sections. Furthermore, {\em only} functions of this form give
rise to such superpositions of singular distributions and of nothing else.

\subsection{Locally Integrable Real Functions}\label{Ssec0203}

In order to deal with the class of non-integrable real functions we are to
discuss here, we must now introduce a different concept of integrability,
which is known as {\em local integrability}. For integrable real functions
in compact domains this is equivalent to both plain integrability and
absolute integrability, since all real functions discussed in this paper
are assumed to be Lebesgue-measurable functions. Here is the definition,
formulated in a way appropriate for our case here.

\begin{definition}\Colon\label{Def01}
  A real function $f(\theta)$ is locally integrable on the unit circle if
  it is integrable on every closed interval contained within that domain.
\end{definition}

\noindent
In this case we say that the real function is {\em locally integrable
  everywhere} on the unit circle. This definition selects the same set of
real functions as plain integrability, and also the same set of real
functions as absolute integrability, given that all real functions
discussed in this paper are assumed to be Lebesgue-measurable functions.
This concept of local integrability can now be generalized to functions
that have a finite number of non-integrable singularities, giving rise to
the concept of {\em local integrability almost everywhere}. Given a real
function $f(\theta)$ defined on the unit circle, which has a finite number
$N\geq 1$ of non-integrable singularities at isolated points on the unit
circle, corresponding to the angles $\{\theta_{1},\ldots,\theta_{N}\}$, we
introduce the definition that follows.

\begin{definition}\Colon\label{Def02}
  A real function $f(\theta)$ is locally integrable almost everywhere on
  the unit circle if it is integrable on every closed interval contained
  within that domain, that does not contain any of the points where the
  function has non-integrable singularities, of which there is a finite
  number.
\end{definition}

\noindent
This definition begins to characterize the class of non-integrable real
functions whose relationship with inner analytic functions we will
consider here. Note that we include in this definition the fact that the
number of non-integrable singularities must be finite.

\subsection{Extended Fourier Theory}\label{Ssec0204}

In~\cite{CAoRFIII} we showed that the whole Fourier theory of integrable
real functions is contained within the complex-analytic structure
presented in~\cite{CAoRFI}. There we also established an extension of the
Fourier theory to essentially all inner analytic functions. In particular,
we introduced the concept of an exponentially bounded sequence of either
Taylor or Fourier coefficients. Here is that definition: given an
arbitrary ordered set of complex coefficients $a_{k}$, for
$k\in\{0,1,2,3,\ldots,\infty\}$, if they satisfy the condition that

\begin{equation}\label{EQExpLim}
  \lim_{k\to\infty}
  |a_{k}|\e{-Ck}
  =
  0,
\end{equation}

\noindent
for all real $C>0$, then we say that the sequence of coefficients $a_{k}$
is {\em exponentially bounded}. This definition applies equally well to
either complex Taylor coefficients or to real Fourier coefficients. In the
construction of the inner analytic functions presented in~\cite{CAoRFI},
if $\alpha_{0}$, $\alpha_{k}$ and $\beta_{k}$, for
$k\in\{1,2,3,\ldots,\infty\}$, are the Fourier coefficients of an
integrable real function $f(\theta)$, then we construct from them the set
of complex Taylor coefficients $c_{0}=\alpha_{0}/2$ and
$c_{k}=\alpha_{k}-\ii\beta_{k}$, for $k\in\{1,2,3,\ldots,\infty\}$. We
also showed in~\cite{CAoRFIII} that the condition that the sequence of
Taylor coefficients $c_{k}$ be exponentially bounded is equivalent to the
condition that the corresponding sequences of Fourier coefficients
$\alpha_{k}$ and $\beta_{k}$ be both exponentially bounded. Another result
that was obtained in~\cite{CAoRFIII} is that the condition that the
sequence of Taylor coefficients $c_{k}$ be exponentially bounded is
sufficient to ensure that the corresponding power series is convergent
within the open unit disk, and therefore converges to an inner analytic
function. Finally, we showed in that paper that the condition that the
Fourier coefficients be exponentially bounded suffices to guarantee that
the {\em regulated} Fourier series given by

\begin{equation}\label{EQRegSer}
  u(\rho,\theta)
  =
  \frac{\alpha_{0}}{2}
  +
  \sum_{k=1}^{\infty}
  \rho^{k}
  \left[
    \alpha_{k}\cos(k\theta)
    +
    \beta_{k}\sin(k\theta)
  \right]
\end{equation}

\noindent
converges absolutely and uniformly for $0<\rho<1$, and that $f(\theta)$
can be obtained as the $\rho\to 1_{(-)}$ limit of this regulated Fourier
series almost everywhere on the unit circle,

\begin{equation}
  f(\theta)
  =
  \frac{\alpha_{0}}{2}
  +
  \lim_{\rho\to 1_{(-)}}
  \sum_{k=1}^{\infty}
  \rho^{k}
  \left[
    \alpha_{k}\cos(k\theta)
    +
    \beta_{k}\sin(k\theta)
  \right].
\end{equation}

\noindent
This has the effect of extending the applicability of the Fourier theory
far beyond the realm of integrable real functions, including, for example,
the singular Schwartz distributions.

\section{General Statement of the Problem}\label{Sec03}

In~\cite{CAoRFI} we showed that every integrable real function can be
related to a unique inner analytic function that corresponds to it. Since
all the inner analytic functions can be organized in integral-differential
chains, it follows that each integrable real function can be assigned a
unique place in a unique integral-differential chain. However, note that
all such functions are assigned only to a part of each one of these
chains, the part where the singularities of the inner analytic functions
are typically either soft or borderline hard. Still, the complete
integral-differential chains do exist and, in particular, do extend
indefinitely in the differentiation direction. As we will see, in general
non-integrable real functions can be related to inner analytic functions
that have one or more hard singularities which are not borderline hard
ones, being therefore on the differentiation side of the
integral-differential chains.

In addition to this, in~\cite{CAoRFII} we showed that integrable real
functions that are piecewise polynomial, and only those, belong to
integral-differential chains that have only superpositions of one or more
delta ``functions'' and derivatives of delta ``functions'' on their
differentiation side, if one goes in the differentiation direction far
enough along their integral-differential chains. Therefore, given any
integrable real function which is {\em not} a piecewise polynomial
function, and is such that the corresponding inner analytic function has
one or more singularities on the unit circle, a type of real function of
which there certainly are many, the integral-differential chain to which
it belongs contains, if we travel far enough in the differentiation
direction along the chain, real functions which are {\em not} integrable
at one or more of their singular points.

Starting from an integrable real function $f(\theta)$ with Fourier
coefficients $\alpha_{k}$ and $\beta_{k}$, leading to the complex Taylor
coefficients $c_{k}$ of an inner analytic function $w(z)$, thus following
the construction presented in~\cite{CAoRFI}, it is not difficult to see
from the definition of the angular derivative that after $n$ angular
differentiations of $w(z)$ we get a proper inner analytic function
$w^{n\ldot}(z)$ with coefficients given by
$c_{k}^{(n)}=\ii^{n}k^{n}c_{k}$. Starting from the convergent Taylor
series of $w(z)$ we have for its $n^{\rm th}$ angular derivative

\noindent
\begin{eqnarray}
  w(z)
  & = &
  \sum_{k=0}^{\infty}
  c_{k}z^{k}
  \;\;\;\Rightarrow
  \nonumber\\
  w^{n\ldot}(z)
  & = &
  \sum_{k=0}^{\infty}
  c_{k}^{(n)}z^{k}
  \nonumber\\
  & = &
  \ii^{n}
  \sum_{k=0}^{\infty}
  \left(k^{n}c_{k}\right)z^{k},
\end{eqnarray}

\noindent
for $n\in\{0,1,2,3,\ldots,\infty\}$. As we showed in~\cite{CAoRFI}, the
coefficients $c_{k}$ associated to an integrable real function $f(\theta)$
are limited as functions of $k$. For integrable functions $f(\theta)$
which are not analytic on the whole unit circle the coefficients $c_{k}$
typically go to zero when $k\to\infty$, at a pace which is slower than
exponentially fast. We thus see that non-integrable real functions
obtained as derivatives of integrable real functions are associated to
Taylor coefficients which typically diverge to infinity as a positive
power of $k$ in the $k\to\infty$ limit.

Therefore, we must conclude that there is a large class of non-integrable
real functions which are still related to corresponding inner analytic
functions such as $w^{n\ldot}(z)$, and that can be obtained as the
$\rho\to 1_{(-)}$ limits to the unit circle of the real parts of these
inner analytic functions. This is the case so long as the $\rho\to
1_{(-)}$ limits of these inner analytic functions exist almost everywhere
on the unit circle. Since angular differentiation conserves the set of
singular points of $w(z)$, at all points where $w(z)$ is analytic, so is
$w^{n\ldot}(z)$, and therefore the limit exists there. Where $w(z)$ has
soft singularities with a large enough degree of softness, namely larger
that or equal to $n+1$, the angular derivative $w^{n\ldot}(z)$ will also
have soft singularities, and again the limit exists there. Since the
$\rho\to 1_{(-)}$ limits do not exist at hard singularities, and in
analogy to what we discussed in~\cite{CAoRFI} in the case of the
borderline hard singularities of integrable real functions, we see that we
must impose some limitations on the non-integrable real functions to be
considered here. Specifically, we assume that the numbers of
non-integrable singularities of these functions are finite.

By the argument above we can immediately determine the existence of many
non-integrable real functions which are related to given inner analytic
functions. Given any inner analytic function which is not simply a
superposition of inner analytic functions associated to singular
distributions, and which has one or more singularities on the unit circle,
by differentiating it a sufficient number of times we produce an inner
analytic function which corresponds to a non-integrable real function.
However, given only a non-integrable real function, the problem of the
determination of the corresponding inner analytic function is not so
immediate. In fact, the problem stated in this way does not really have a
unique solution, because such a real function is not really defined on the
whole periodic interval, but only on a strict subset of it. The definition
of a non-integrable real function must leave out the set of points on the
unit circle where it has non-integrable singularities, and thus diverges
to infinity. Therefore, one can superpose to such a real function any set
of delta ``functions'' and derivatives of delta ``functions'', that have
their singularities located at those same singular points, without
changing at all the definition of the original non-integrable real
function. It then follows that we can add to any inner analytic function
that corresponds to the non-integrable real function linear combinations
of the inner analytic functions that correspond to the delta ``function''
and to the derivatives of the delta ``function'', which were given
explicitly in~\cite{CAoRFII}, chosen so that their singularities are
located at the points where the non-integrable real function is not
defined, without changing the correspondence between the non-integrable
real function and the respective inner analytic function.

Let us state, then, the problem we propose to investigate in this paper,
relating to non-integrable real functions and corresponding inner analytic
functions. We want to find out how to determine an inner analytic
function, given only the non-integrable real function, that reproduces
that real function almost everywhere in its domain of definition as the
$\rho\to 1_{(-)}$ limit of its real part. We will give here a solution to
this problem for a rather large class of non-integrable real functions, on
which we will impose, however, some limitations. The first limitation will
be, of course, that the real function $f(\theta)$ be locally integrable
almost everywhere. This means, in particular, that its number of hard
singular points, either integrable or non-integrable, must be finite. The
other condition is that the non-integrable singular points of $f(\theta)$
have finite degrees of hardness, as defined and discussed for inner
analytic functions in~\cite{CAoRFI}, and translated to real functions in
this paper, at the beginning of Section~\ref{Sec02}.

\section{Representation of Non-Integrable Real Functions}\label{Sec04}

First of all, let us describe, in general lines, the algorithm we propose
to use in order to determine an inner analytic function that corresponds
to a given non-integrable real function. Given a non-integrable real
function $f(\theta)$, which is however locally integrable almost
everywhere, and whose hard singularities have a finite maximum degree of
hardness $n$, we sectionally integrate it $n$ times. Since by hypothesis
all the non-integrable singularities of $f(\theta)$ have degrees of
hardness of at most $n$, the resulting function $f^{-n\prime}(\theta)$ is
in fact an integrable one on the whole unit circle. Therefore, we may use
it to construct the corresponding inner analytic function, which we will
name $w^{-n\ldot}(z)$, using the construction presented in~\cite{CAoRFI}.
Having this inner analytic function, we then calculate its $n^{\rm th}$
angular derivative, in order to obtain $w(z)$, which is the inner analytic
function that corresponds to the non-integrable real function $f(\theta)$.
However, the $n$-fold angular differentiation process produces only the
{\em proper} inner analytic function $w_{p}(z)$ associated to $w(z)$, and
therefore produces $f(\theta)$ only up to an overall constant related to
the whole unit circle, as we will soon see. Of course this scheme will
succeed if and only if the non-integrable singular points of $f(\theta)$
all have finite degrees of hardness, with a maximum value of $n$.

\vspace{2.6ex}

\noindent
In this section we will prove the following theorem.

\begin{theorem}\Colon\label{Theo01}
  Every non-integrable real function defined almost everywhere on the
  periodic interval, which is locally integrable almost everywhere, and
  which is such that its non-integrable singularities have finite degrees
  of hardness, can be represented by an inner analytic function, and can
  be recovered almost everywhere on its domain of definition by means of
  the limit to the unit circle of the real part of that inner analytic
  function.
\end{theorem}

\vspace{2.6ex}

\noindent
Before we attempt to prove the theorem, let us establish some notation for
a sectionally defined real function $f(\theta)$, that will be similar to
the one adopted for the piecewise polynomial functions discussed
in~\cite{CAoRFII}, which is also given in Equation~(\ref{EQPieFun}).
Since a non-integrable real function $f(\theta)$ which is locally
integrable almost everywhere is not defined at the singular points
corresponding to the angles $\theta_{i}$, for $i\in\{1,\ldots,N\}$, it is
in fact only sectionally defined, in $N$ open intervals between
consecutive singularities, as given in Equation~(\ref{EQSeqSec}). Let us
therefore specify the definition of $f(\theta)$ as a set of $N\geq 1$
functions $f_{i}(\theta)$ defined on the $N$ sections
$(\theta_{i},\theta_{i+1})$,

\begin{equation}
  f(\theta)
  =
  \left\{
    \rule{0em}{2.5ex}
    f_{i}(\theta),i\in\{1,\ldots,N\}
  \right\},
\end{equation}

\noindent
where, as before, we adopt the convention that every section
$(\theta_{i},\theta_{i+1})$ is numbered after the singular point
$\theta_{i}$ at its left end.

As a preliminary to the proof of the theorem, some considerations are in
order, regarding the multiple sectional integration of such non-integrable
functions. Note that, if $f(\theta)$ is locally integrable almost
everywhere, then it is integrable within each one of the $N$ open
intervals that define the sections. More precisely, it is integrable on
every closed interval contained within one of these open intervals. In
terms of the sectional functions, since $f_{i}(\theta)$ is integrable
within its section, we may define a piecewise primitive for $f(\theta)$,
by simply integrating each function $f_{i}(\theta)$ within the
corresponding open interval, starting at some arbitrary reference point
$\theta_{0,i}$ strictly within that open interval. This will define the
piecewise primitive

\noindent
\begin{eqnarray}
  f^{-1\prime}(\theta)
  & = &
  \left\{
    \rule{0em}{2.5ex}
    f_{i}^{-1\prime}(\theta),i\in\{1,\ldots,N\}
  \right\},
  \nonumber\\
  f_{i}^{-1\prime}(\theta)
  & = &
  \int_{\theta_{0,i}}^{\theta}d\theta'\,
  f_{i}(\theta'),
\end{eqnarray}

\noindent
where we have that $\theta_{i}<\theta<\theta_{i+1}$ and that
$\theta_{i}<\theta_{0,i}<\theta_{i+1}$. Note that, since $f_{i}(\theta)$
is integrable on every closed interval contained within its section, it
follows that $f_{i}^{-1\prime}(\theta)$ is limited on these closed
intervals, and therefore is also integrable on them. Therefore, this
process of sectional integration of the sectional functions can be
iterated indefinitely. If we iterate this process of sectional
integration, we obtain further piecewise primitives
$f^{-2\prime}(\theta)$, $f^{-3\prime}(\theta)$, and so on. Note also that,
during this process of multiple sectional integration, some non-integrable
singularities, having a lower degree of hardness, may become integrable
before the others. In this case we might ignore the singularities which
became integrable, from that point on in the multiple integration process,
which therefore effectively reduces the number of sections, but for
simplicity we choose not to do that, and thus to keep the set of sections
constant. Note that, in any case, we do continue the process of sectional
integration until {\em all} the singularities have become integrable, of
course.

Although for definiteness we are integrating from some particular
reference points $\theta_{0,i}$ within each section, our objective here is
actually to construct primitives, and therefore we may ignore the
particular values chosen for $\theta_{0,i}$ if at each step in this
iterative process we add an arbitrary integration constant to the
primitive in each section, so that after $n$ such successive integrations
a polynomial of degree $n-1$, with $n$ arbitrary coefficients, will have
been added to the $n^{\rm th}$ primitive in the $i^{\rm th}$ section. We
will express this as follows,

\noindent
\begin{eqnarray}
  f^{-n\prime}(\theta)
  & = &
  \left\{
    \rule{0em}{2.5ex}
    f_{i}^{-n\prime}(\theta)
    +
    P_{i}^{(n-1)}(\theta),
    i\in\{1,\ldots,N\}
  \right\}
  \nonumber\\
  & = &
  \left\{
    \rule{0em}{2.5ex}
    f_{i}^{-n\prime}(\theta),
    i\in\{1,\ldots,N\}
  \right\}
  +
  P_{(n-1)}(\theta),
  \nonumber\\
  P_{(n-1)}(\theta)
  & = &
  \left\{
    \rule{0em}{2.5ex}
    P_{i}^{(n-1)}(\theta),
    i\in\{1,\ldots,N\}
  \right\},
\end{eqnarray}

\noindent
where $f^{-n\prime}(\theta)$ is the most general piecewise $n^{\rm th}$
primitive of $f(\theta)$, $f_{i}^{-n\prime}(\theta)$ is an arbitrary
$n^{\rm th}$ primitive of $f_{i}(\theta)$ in the $i^{\rm th}$ section,
$P_{i}^{(n-1)}(\theta)$ is an arbitrary polynomial of order $n-1$ in the
$i^{\rm th}$ section, and $P_{(n-1)}(\theta)$ is a piecewise polynomial
function of order $n-1$, containing therefore $n$ arbitrary constants in
each section. Note that $P_{(n-1)}(\theta)$ is always an integrable real
function. Note also that, upon subsequent $n$-fold differentiation of
$f^{-n\prime}(\theta)$ with respect to $\theta$, all the arbitrary
constants that were added during the multiple integration process are then
eliminated, the arbitrary polynomials vanish from the result, and we thus
get back the original function $f(\theta)$, within the open intervals that
constitute the sections,

Finally, we must emphasize some facts about the behavior of the
correspondence between inner analytic functions an integrable real
functions under the respective operations of differentiation and
integration, that take us along the corresponding integral-differential
chain. First, let us recall that, as we saw in~\cite{CAoRFI}, both angular
differentiation and angular integration produce only {\em proper} inner
analytic functions, and thus always result in null Taylor coefficients
$c_{0}$, and thus in null Fourier coefficients $\alpha_{0}$. Therefore,
when we use angular integration and differentiation in our algorithm, we
lose all information about these two $k=0$ coefficients. Second, as we
also saw in~\cite{CAoRFI}, the integration on $\theta$ implies that the
resulting Fourier coefficient $\alpha_{0}$ becomes indeterminate due to
the presence of an arbitrary integration constant. It is due to this that
we must add arbitrary constants during the process of multiple sectional
integration, thus generating the arbitrary piecewise polynomial real
function $P_{(n-1)}(\theta)$.

At this point, let us review the algorithm we are to use here. First,
starting from the non-integrable real function $f(\theta)$ we go $n$ steps
along the integration direction of the integral-differential chain, using
sectional integration on $\theta$ on the unit circle. This produces the
$n$-fold piecewise primitive $f^{-n\prime}(\theta)$ containing the
arbitrary piecewise polynomial real function $P_{(n-1)}(\theta)$. From the
globally integrable real function $f^{-n\prime}(\theta)$ we then construct
the inner analytic function $w^{-n\ldot}(z)$, using the construction
presented in~\cite{CAoRFI}. We then come back in the differentiation
direction of the integral-differential chain the same number of steps,
using this time angular differentiation of the inner analytic functions
within the open unit disk. This produces an inner analytic function
associated to $f(\theta)$, except for its coefficient $c_{0}$, which means
that we recover only a proper inner analytic function, given that it is
the result of a series of angular differentiations, and we will therefore
denote this function by $w_{p}(z)=u_{p}(\rho,\theta)+\ii
v_{p}(\rho,\theta)$. Since this is equivalent to differentiation with
respect to $\theta$ of the piecewise polynomial real function
$P_{(n-1)}(\theta)$, it completely eliminates this function within the
open intervals that constitute the sections. However, it also produces the
superposition of a set of delta ``functions'' and derivatives of delta
``functions'' with support on the singular points between successive
sections.

Thus we must conclude that this algorithm necessarily results in an
indeterminate Taylor coefficient $c_{0}$ of the inner analytic function
$w(z)$ which corresponds to the real function $f(\theta)$, so that we
recover only the corresponding proper inner analytic function $w_{p}(z)$,
and therefore it results in an equally indeterminate Fourier coefficients
$\alpha_{0}$ related to the real function $f(\theta)$. Se see, therefore,
that in this paper we are compelled to relax, to some extent, and only
during the process of construction of the inner analytic functions, the
correspondence between the real functions and these inner analytic
functions, considering only proper inner analytic functions, and accepting
the fact that they will correspond to the non-integrable real functions
only up to a global constant over the whole unit circle. Once the proper
inner analytic function $w_{p}(z)$ corresponding to a non-integrable real
function $f(\theta)$ has been determined, the relation will be written as

\begin{equation}
  f(\theta)
  \longleftrightarrow
  \frac{\alpha_{0}}{2}
  +
  w_{p}(z),
\end{equation}

\noindent
where $\alpha_{0}$ is a real constant, which can be determined afterwards,
by the simple comparison of the known value of the real part
$u_{p}(1,\theta)$ of $w_{p}(z)$ and the known value of $f(\theta)$, at any
point $\theta$ on the unit circle where they do not have hard
singularities. If $\theta_{0}$ is such a point, then we have that
$\alpha_{0}=2[f(\theta_{0})-u_{p}(1,\theta_{0})]$. Once the coefficient
$\alpha_{0}$ is thus determined, this finally determines completely the
inner analytic function $w(z)$ that corresponds to $f(\theta)$,

\begin{equation}
  w(z)
  =
  \frac{\alpha_{0}}{2}
  +
  w_{p}(z).
\end{equation}

\noindent
We are now ready to prove the theorem.

\begin{proof}\Colon
\end{proof}

\noindent
In order to prove the theorem, our first task is to show that we can
obtain an inner analytic function $w(z)$ from the non-integrable real
function $f(\theta)$, which is assumed to be locally integrable almost
everywhere. Consider that we execute the iterative piecewise integration
process described before $n$ times on $f(\theta)$, where $n$ is the
maximum among all the degrees of hardness of the non-integrable hard
singularities involved. By doing this we generate a real function
$f^{-n\prime}(\theta)$ which has only soft or at most borderline hard
singularities on the unit circle, and which is therefore integrable on the
whole unit circle.

It follows therefore that we may determine its set of Fourier
coefficients, as was done in~\cite{CAoRFI}, which we will name
$\alpha_{0}^{(-n)}$, $\alpha_{k}^{(-n)}$ and $\beta_{k}^{(-n)}$, for
$k\in\{1,2,3,\ldots,\infty\}$. From this set of Fourier coefficients we
may then define, again as was done in~\cite{CAoRFI}, the complex Taylor
coefficients $c_{k}^{(-n)}$, for $k\in\{0,1,2,3,\ldots,\infty\}$. From
these coefficients we may then determine the unique inner analytic
function that corresponds to the $n^{\rm th}$ piecewise primitive
$f^{-n\prime}(\theta)$, which we will name
$w^{-n\ldot}(z)=u^{-n\prime}(\rho,\theta)+\ii v^{-n\prime}(\rho,\theta)$.
Note that we may use this notation unequivocally because every proper
inner analytic function belongs to an infinite integral-differential
chain, extending indefinitely to either side, so that we know that there
exists in fact a proper inner analytic function $w_{p}(z)$ associated to

\begin{equation}
  w_{p}^{-n\ldot}(z)
  =
  w^{-n\ldot}(z)
  -
  c_{0}^{(-n)},
\end{equation}

\noindent
namely its $n^{\rm th}$ angular derivative. As we have shown
in~\cite{CAoRFI}, from the $\rho\to 1_{(-)}$ limit to the unit circle of
the real part $u^{-n\prime}(\rho,\theta)$ of the inner analytic function
$w^{-n\ldot}(z)$ we can recover the integrable real function
$f^{-n\prime}(\theta)$, almost everywhere in its domain of definition. We
thus have the correspondence for the integrable real function

\begin{equation}
  f^{-n\prime}(\theta)
  \longleftrightarrow
  w^{-n\ldot}(z).
\end{equation}

\noindent
Having done this, we now take the $n^{\rm th}$ angular derivative of the
inner analytic function $w^{-n\ldot}(z)$, thus obtaining a proper inner
analytic function $w_{p}(z)=u_{p}(\rho,\theta)+\ii
v_{p}(\rho,\theta)$. Since $n$ angular differentiations of
$w^{-n\ldot}(z)$ correspond to $n$ differentiations with respect to
$\theta$ of $f^{-n\prime}(\theta)$, and thus completely eliminates the
piecewise polynomial real function $P_{(n-1)}(\theta)$ within the domain
of definition of $f(\theta)$, it follows that the function $w_{p}(z)$ is a
proper inner analytic function corresponding to the non-integrable real
function $f(\theta)$. As discussed before, this correspondence is valid
only up to a global real constant yet to be determined, so that we may now
write

\begin{equation}
  f(\theta)
  \longrightarrow
  \frac{\alpha_{0}}{2}
  +
  w_{p}(z),
\end{equation}

\noindent
where $\alpha_{0}$ can then be determined as discussed before, by the
comparison between the real part $u_{p}(1,\theta)$ of $w_{p}(z)$ and
$f(\theta)$ at some particular point on the unit circle where they do not
have hard singularities. Having determined $\alpha_{0}$, and therefore
$c_{0}=\alpha_{0}/2$, we may now define the inner analytic function that
corresponds to $f(\theta)$,

\begin{equation}
  w(z)
  =
  c_{0}
  +
  w_{p}(z),
\end{equation}

\noindent
so that we have the relation leading from $f(\theta)$ to $w(z)$

\begin{equation}
  f(\theta)
  \longrightarrow
  w(z).
\end{equation}

\noindent
This completes the first part of the proof of Theorem~\ref{Theo01}.

\begin{proof}\Colon
\end{proof}

\noindent
We must now prove that we can recover $f(\theta)$ from the real part
$u(\rho,\theta)$ of $w(z)$ in the $\rho\to 1_{(-)}$ limit. In order to do
this, we start from the fact that from~\cite{CAoRFI} we know this to be
true for the $n$-fold primitives

\begin{equation}
  w^{-n\ldot}(z)
  \longleftrightarrow
  f^{-n\prime}(\theta).
\end{equation}

\noindent
While the inner analytic function $w^{-n\ldot}(z)$ is given by the power
series

\begin{equation}
  w^{-n\ldot}(z)
  =
  \sum_{k=0}^{\infty}
  c_{k}^{(-n)}z^{k},
\end{equation}

\noindent
as we have shown in~\cite{CAoRFIII} the real function
$f^{-n\prime}(\theta)$ can be expressed almost everywhere as an regulated
Fourier series, even if the Fourier series itself diverges almost
everywhere,

\begin{equation}
  f^{-n\prime}(\theta)
  =
  \frac{\alpha_{0}^{(-n)}}{2}
  +
  \lim_{\rho\to 1_{(-)}}
  \sum_{k=1}^{\infty}
  \rho^{k}
  \left[
    \alpha_{k}^{(-n)}\cos(k\theta)
    +
    \beta_{k}^{(-n)}\sin(k\theta)
  \right],
\end{equation}

\noindent
where we have that $c_{0}^{(-n)}=\alpha_{0}^{(-n)}/2$ and that
$c_{k}^{(-n)}=\alpha_{k}^{(-n)}-\ii\beta_{k}^{(-n)}$ for
$k\in\{1,2,3,\ldots,\infty\}$. As we established in~\cite{CAoRFIII}, since
the sequences of Fourier coefficients $\alpha_{k}^{(-n)}$ and
$\beta_{k}^{(-n)}$ are exponentially bounded, this series is absolutely
and uniformly convergent for $0<\rho<1$. As we saw in~\cite{CAoRFI}, the
fact that we can recover $f^{-n\prime}(\theta)$ as the $\rho\to 1_{(-)}$
limit of the real part $u^{-n\prime}(\rho,\theta)$ of $w^{-n\ldot}(z)$ is
a consequence of the fact that the two real functions
$u^{-n\prime}(1,\theta)$ and $f^{-n\prime}(\theta)$ have exactly the same
set of Fourier coefficients. This, is turn, can be expressed as the
relations between the Taylor coefficients $c_{k}^{(-n)}$ associated to
$u(\rho,\theta)$ and the Fourier coefficients $\alpha_{k}^{(-n)}$ and
$\beta_{k}^{(-n)}$ associated to $f(\theta)$, which have just been given
above.

Let us now prove that the correspondence between
$u^{-n\prime}(\rho,\theta)$ and $f^{-n\prime}(\theta)$ implies the same
correspondence between $u^{(-n+1)\prime}(\rho,\theta)$ and
$f^{(-n+1)\prime}(\theta)$, when we differentiate the two functions. We
can do this by just showing that the relations between the Taylor
coefficients and the Fourier coefficients are preserved by this process of
differentiation. If we just differentiate the inner analytic function
using angular differentiation we get

\noindent
\begin{eqnarray}
  w^{(-n+1)\ldot}(z)
  & = &
  \sum_{k=1}^{\infty}
  \ii k\,c_{k}^{(-n)}z^{k}
  \nonumber\\
  & = &
  \sum_{k=1}^{\infty}
  c_{k}^{(-n+1)}z^{k},
\end{eqnarray}

\noindent
from which we have for the Taylor coefficients $c_{k}^{(-n+1)}=\ii
k\,c_{k}^{(-n)}$, for $k\in\{1,2,3,\ldots,\infty\}$. Note that, since this
is a convergent power series, we can always differentiate it
term-by-term. If we now differentiate the real function using simple
differentiation with respect to $\theta$ we get

\noindent
\begin{eqnarray}
  f^{(-n+1)\prime}(\theta)
  & = &
  \lim_{\rho\to 1_{(-)}}
  \sum_{k=1}^{\infty}
  \rho^{k}
  \left[
    k\,\beta_{k}^{(-n)}\cos(k\theta)
    -
    k\,\alpha_{k}^{(-n)}\sin(k\theta)
  \right],
  \nonumber\\
  & = &
  \lim_{\rho\to 1_{(-)}}
  \sum_{k=1}^{\infty}
  \rho^{k}
  \left[
    \alpha_{k}^{(-n+1)}\cos(k\theta)
    +
    \beta_{k}^{(-n+1)}\sin(k\theta)
  \right],
\end{eqnarray}

\noindent
from which we have for the corresponding Fourier coefficients that
$\alpha_{k}^{(-n+1)}=k\,\beta_{k}^{(-n)}$ and also that
$\beta_{k}^{(-n+1)}=-k\,\alpha_{k}^{(-n)}$, for
$k\in\{1,2,3,\ldots,\infty\}$. Note that we have, in either case, that
$c_{0}^{(-n+1)}=0$ and that $\alpha_{0}^{(-n+1)}=0$, so that the relation
between $c_{0}^{(-n)}$ and $\alpha_{0}^{(-n)}$ is in fact preserved. Note
also that, since this trigonometric series is uniformly convergent, and in
fact is the real part of a convergent complex power series, we may
differentiate it term-by-term so long as the series thus obtained is also
convergent. Since the Fourier coefficients $\alpha_{0}^{(-n+1)}$ and
$\beta_{0}^{(-n+1)}$, which increase at most as a power of $k$ when
$k\to\infty$, are thus seen to be exponentially bounded, this implies that
the series thus obtained is also absolutely and uniformly convergent, as
we have shown in~\cite{CAoRFIII}. Therefore, we are justified in
differentiating the original series term-by-term. We therefore have the
relation between the Taylor coefficients $c_{k}^{(-n+1)}$ and the Fourier
coefficients $\alpha_{k}^{(-n+1)}$ and $\beta_{k}^{(-n+1)}$,

\noindent
\begin{eqnarray}
  c_{k}^{(-n+1)}
  & = &
  \ii k\,c_{k}^{(-n)}
  \nonumber\\
  & = &
  \ii k
  \left[
    \alpha_{k}^{(-n)}
    -
    \ii
    \beta_{k}^{(-n)}
  \right]
  \nonumber\\
  & = &
  \ii k
  \left[
    -\,
    \frac{1}{k}\,
    \beta_{k}^{(-n+1)}
    -
    \ii\,
    \frac{1}{k}\,
    \alpha_{k}^{(-n+1)}
  \right]
  \nonumber\\
  & = &
  \alpha_{k}^{(-n+1)}
  -
  \ii
  \beta_{k}^{(-n+1)}.
\end{eqnarray}

\noindent
We see therefore that we indeed have that the relation between
$c_{k}^{(-n)}$, $\alpha_{k}^{(-n)}$ and $\beta_{k}^{(-n)}$, for
$k\in\{1,2,3,\ldots,\infty\}$, is also preserved,

\begin{equation}
  c_{k}^{(-n+1)}
  =
  \alpha_{k}^{(-n+1)}
  -
  \ii
  \beta_{k}^{(-n+1)},
\end{equation}

\noindent
which thus establishes the correspondence for the first derivatives. We
may now extend this argument to subsequent derivatives, from
$w^{(-n+1)\ldot}(z)$ and $f^{(-n+1)\prime}(\theta)$ all the way to $w(z)$
and $f(\theta)$, by finite induction. Therefore, let us assume the result
for the case $(-n+i)$ and show that this implies that it is also valid for
the case $(-n+i+1)$. We assume therefore that we have

\begin{equation}
  c_{k}^{(-n+i)}
  =
  \alpha_{k}^{(-n+i)}
  -
  \ii
  \beta_{k}^{(-n+i)},
\end{equation}

\noindent
for some positive value of $i$, where the functions are expressed as the
corresponding series

\noindent
\begin{eqnarray}
  w^{(-n+i)\ldot}(z)
  & = &
  \sum_{k=1}^{\infty}
  c_{k}^{(-n+i)}z^{k},
  \nonumber\\
  f^{(-n+i)\prime}(\theta)
  & = &
  \lim_{\rho\to 1_{(-)}}
  \sum_{k=1}^{\infty}
  \rho^{k}
  \left[
    \alpha_{k}^{(-n+i)}\cos(k\theta)
    +
    \beta_{k}^{(-n+i)}\sin(k\theta)
  \right].
\end{eqnarray}

\noindent
We may now differentiate either series term-by-term, which we may do for
the same reasons as before, thus obtaining

\noindent
\begin{eqnarray}
  w^{(-n+i+1)\ldot}(z)
  & = &
  \sum_{k=1}^{\infty}
  \ii k\,c_{k}^{(-n+i)}z^{k},
  \nonumber\\
  f^{(-n+i+1)\prime}(\theta)
  & = &
  \lim_{\rho\to 1_{(-)}}
  \sum_{k=1}^{\infty}
  \rho^{k}
  \left[
    k\,\beta_{k}^{(-n+i)}\cos(k\theta)
    -
    k\,\alpha_{k}^{(-n+i)}\sin(k\theta)
  \right],
\end{eqnarray}

\noindent
so that we have for the coefficients for the case $(-n+i+1)$

\noindent
\begin{eqnarray}
  c_{k}^{(-n+i+1)}
  & = &
  \ii k\,c_{k}^{(-n+i)},
  \nonumber\\
  \alpha_{k}^{(-n+i+1)}
  & = &
  k\,\beta_{k}^{(-n+i)},
  \nonumber\\
  \beta_{k}^{(-n+i+1)}
  & = &
  -
  k\,\alpha_{k}^{(-n+i)},
\end{eqnarray}

\noindent
for $k\in\{1,2,3,\ldots,\infty\}$. Using now the relations between the
coefficients for the case $(-n+i)$ we have

\noindent
\begin{eqnarray}
  c_{k}^{(-n+i+1)}
  & = &
  \ii k\,c_{k}^{(-n+i)}
  \nonumber\\
  & = &
  \ii k
  \left[
    \alpha_{k}^{(-n+i)}
    -
    \ii
    \beta_{k}^{(-n+i)}
  \right]
  \nonumber\\
  & = &
  \ii k
  \left[
    -\,
    \frac{1}{k}\,
    \beta_{k}^{(-n+i+1)}
    -
    \ii\,
    \frac{1}{k}\,
    \alpha_{k}^{(-n+i+1)}
  \right]
  \nonumber\\
  & = &
  \alpha_{k}^{(-n+i+1)}
  -
  \ii
  \beta_{k}^{(-n+i+1)},
\end{eqnarray}

\noindent
for $k\in\{1,2,3,\ldots,\infty\}$, thus showing that the relation between
the coefficients is in fact preserved, where we recall that the $k=0$
coefficients are always zero during this process. This is therefore true
for all possible multiple derivatives, all the way to infinity, and in
particular it is true for $i=n$, that is, for the coefficients of $w(z)$
and $f(\theta)$. In order to complete the proof in this case all we have
to do is to consider the real function

\begin{equation}
  g(\theta)
  =
  u(1,\theta)-f(\theta),
\end{equation}

\noindent
where

\begin{equation}
  u(1,\theta)
  =
  \lim_{\rho\to 1_{(-)}}
  u(\rho,\theta).
\end{equation}

\noindent
Since the expression of the Fourier coefficients is linear on the
functions, and since $u(1,\theta)$ and $f(\theta)$ have exactly the same
set of Fourier coefficients, it is clear that all the Fourier coefficients
of $g(\theta)$ are zero. Therefore, for the real function $g(\theta)$ we
have that $c_{k}=0$ for all $k$, and thus the inner analytic function that
corresponds to $g(\theta)$ is the identically null complex function
$w_{\gamma}(z)\equiv 0$. This is an inner analytic function which is, in
fact, analytic over the whole complex plane, and which, in particular, is
zero over the unit circle, so that we have $g(\theta)\equiv 0$, since the
identically zero real function is the {\em only} integrable real function
without removable singularities that corresponds to the identically zero
inner analytic function, due to the completeness of the Fourier basis, as
was shown in~\cite{CAoRFIII}. Note, in particular, that the $\rho\to
1_{(-)}$ limits exist at {\em all} points of the unit circle in the case
of the inner analytic function associated to $g(\theta)$. Since both
$u(1,\theta)$ and $f(\theta)$ have non-integrable hard singularities at
isolated points on the unit circle, we can conclude only that

\begin{equation}
  f(\theta)
  =
  \lim_{\rho\to 1_{(-)}}
  u(\rho,\theta)
\end{equation}

\noindent
almost everywhere on the unit circle, or everywhere in the domain of
definition of $f(\theta)$, which therefore excludes all the points where
the function has hard singularities. Note that the domain of definition of
$u(1,\theta)$ is the same as that of $f(\theta)$, because any hard
singularities that may have been softened in the process of iterative
integration, during the construction of $w(z)$, will have been hardened
again in the corresponding process of iterative differentiation. We have
therefore the complete correspondence

\begin{equation}
  f(\theta)
  \longleftrightarrow
  w(z).
\end{equation}

\noindent
The inner analytic function $w(z)$ represents the non-integrable real
function $f(\theta)$ exactly in the same way as that which was established
for integrable real functions in~\cite{CAoRFI}. Note that this proof
establishes that the correspondence between the inner analytic functions
and the real functions is also valid for all the intermediate cases, from
$w^{-n\ldot}(z)$ and $f^{-n\prime}(\theta)$ to $w(z)$ and $f(\theta)$.
This completes the second part of the proof of Theorem~\ref{Theo01}.

\vspace{2.6ex}

\noindent
Let us once more draw attention to the fact that the inner analytic
function $w(z)$ produced in the way described above, while it does
correspond to the non-integrable real function $f(\theta)$, is {\em not}
unique. One can add to it any linear combination of inner analytic
functions that correspond to singular distributions with their
singularities at the singular points on the unit circle corresponding to
the angles $\{\theta_{1},\ldots,\theta_{N}\}$, without changing the fact
that the non-integrable real function $f(\theta)$ is still recovered as
the $\rho\to 1_{(-)}$ limit of the real part of the new inner analytic
function obtained in this way, at all points where it is well defined.

Given the inner analytic function $w(z)$ that was obtained from
$f(\theta)$ by the process described above, one can then consider defining
from it a {\em reduced} inner analytic function $w_{r}(z)$ that represents
the non-integrable real function in a unique way, without the
superposition of any singular distributions. If the $N$ singular points of
$w(z)$ are examined in order to determine the existence there of singular
distributions, and since each one of these singular distributions is
represented by a known unique inner analytic function of its own, as was
shown in~\cite{CAoRFII}, one can simply subtract from $w(z)$ the
appropriate linear combination of inner analytic functions of the singular
distributions present, in order to obtain an inner analytic function that
represents the non-integrable real function $f(\theta)$ alone, without any
singular distributions superposed to it.

The simplest way to do this is to examine the result of every angular
differentiation during the process leading from $w^{-n\ldot}(z)$ to
$w(z)$. At each step one can verify whether or not the angular
differentiation has generated one or more Dirac delta ``functions'' at
some of the singular points. This is a simple thing to verify, because the
occurrence of a delta ``function'' at a certain point is always preceded
by the occurrence of a finite discontinuity of the real function at that
point. This can also be done by the determination of the type and
orientation of the singularities at these points, since the Dirac delta
``functions'' are associated to inner analytic functions with simple poles
that have a specific orientation with respect to the unit circle. One can
then subtract from the proper inner analytic function obtained at that
point in the iterative differentiation process the inner analytic
functions corresponding to the delta ``functions'' at the appropriate
points. By doing this one guarantees that no derivatives of delta
``functions'' will ever arise during the process of multiple
differentiation. After one determines $\alpha_{0}$ at the last step of the
process, this will lead to a reduced inner analytic function $w_{r}(z)$
which includes no singular distributions at all, and that hence
corresponds to $f(\theta)$ in a unique and simple way,

\begin{equation}
  f(\theta)
  \longleftrightarrow
  w_{r}(z),
\end{equation}

\noindent
that does not include the superposition of any singular distributions with
support on the singular points between successive sections.

\section{Representation in the Extended Fourier Theory}\label{Sec05}

Assuming that, given a non-integrable real function $f(\theta)$, the
corresponding reduced inner analytic function
$w_{r}(z)=u_{r}(\rho,\theta)+\ii v_{r}(\rho,\theta)$ has been determined,
we may now consider determining a unique set of Fourier coefficients to be
associated to the non-integrable real function $f(\theta)$. Of these,
$\alpha_{0}$ has already been determined, through the comparison of
$u_{r}(1,\theta)$ and $f(\theta)$ at some point of the unit circle where
they do not have hard singularities. It follows that $c_{0}=\alpha_{0}/2$
has also been determined. From the Taylor series of $w_{r}(z)$,

\begin{equation}
  w_{r}(z)
  =
  \sum_{k=0}^{\infty}
  c_{k}z^{k},
\end{equation}

\noindent
we have the values of all the other Taylor coefficients $c_{k}$, for
$k\in\{1,2,3,\ldots,\infty\}$. Given that in this case we have that
$c_{k}=\alpha_{k}-\ii\beta_{k}$, we immediately get the values of all the
other Fourier coefficients $\alpha_{k}$ and $\beta_{k}$. Note that this
construction has the effect of associating a complete set of Fourier
coefficients $\alpha_{0}$, $\alpha_{k}$ and $\beta_{k}$, for
$k\in\{1,2,3,\ldots,\infty\}$, to the non-integrable real function
$f(\theta)$. In particular, the association of $\alpha_{0}$ has the effect
of attributing an {\em average value} to the non-integrable real function
$f(\theta)$, which has been defined via an analytic process.

Although these coefficients obviously cannot be written in terms of
integrals involving $f(\theta)$ in the usual way, all of them except for
$\alpha_{0}$ can, in fact, be written as a certain set of integrals. In
order to do this, we start by considering the $n^{\rm th}$ angular
primitive $w_{r}^{-n\ldot}(z)$ of the reduced inner analytic function
$w_{r}(z)$, which is, therefore, a proper inner analytic function, and the
known associated Fourier and Taylor coefficients, which we will name
$\alpha_{0}^{(-n)}$, $\alpha_{k}^{(-n)}$ and $\beta_{k}^{(-n)}$, for
$k\in\{1,2,3,\ldots,\infty\}$, and $c_{k}^{(-n)}$, for
$k\in\{0,1,2,3,\ldots,\infty\}$. Since it has at most borderline hard
singularities on the unit circle, this inner analytic function corresponds
to an integrable real function $f_{r}^{-n\prime}(\theta)$ on that circle,

\begin{equation}
  f_{r}^{-n\prime}(\theta)
  =
  \lim_{\rho\to 1_{(-)}}
  u_{r}^{-n\ldot}(\rho,\theta),
\end{equation}

\noindent
so that the Fourier coefficients of $f_{r}^{-n\prime}(\theta)$ can be
written as integrals involving this function,

\noindent
\begin{eqnarray}
  \alpha_{k}^{(-n)}
  & = &
  \frac{1}{\pi}
  \int_{-\pi}^{\pi}d\theta\,
  \cos(k\theta)f_{r}^{-n\prime}(\theta),
  \nonumber\\
  \beta_{k}^{(-n)}
  & = &
  \frac{1}{\pi}
  \int_{-\pi}^{\pi}d\theta\,
  \sin(k\theta)f_{r}^{-n\prime}(\theta),
\end{eqnarray}

\noindent
for $k\in\{1,2,3,\ldots,\infty\}$. We now recall that we have for the
Taylor series of the $n^{\rm th}$ angular primitive of $w_{r}(z)$,

\begin{equation}
  w_{r}^{-n\ldot}(z)
  =
  \sum_{k=1}^{\infty}
  c_{k}^{(-n)}z^{k},
\end{equation}

\noindent
and therefore the $n^{\rm th}$ angular derivative of this equation is
given by

\begin{equation}
  w_{r}^{0\ldot}(z)
  =
  \sum_{k=1}^{\infty}
  \left[
    \ii^{n}k^{n}c_{k}^{(-n)}
  \right]z^{k},
\end{equation}

\noindent
where $w_{r}^{0\ldot}(z)$ is the proper inner analytic function associated
to $w_{r}(z)$. It thus follows that we have for the Taylor coefficients
$c_{k}$ of $w_{r}^{0\ldot}(z)$, which are also the Taylor coefficients of
$w_{r}(z)$, for $k\in\{1,2,3,\ldots,\infty\}$,

\begin{equation}\label{EQGenCoe}
  c_{k}
  =
  \ii^{n}k^{n}c_{k}^{(-n)}.
\end{equation}

\noindent
Since we know that the coefficients $c_{k}^{(-n)}$, being the Taylor
coefficients associated to an integrable real function, are limited for
all $k$, we may now conclude that the coefficients $c_{k}$ of the
non-integrable real function may diverge with $k$, but not faster than the
power $k^{n}$. We may write this relation in terms of the Fourier
coefficients $\alpha_{k}$ and $\beta_{k}$ associated to the non-integrable
real function $f(\theta)$, for $k\in\{1,2,3,\ldots,\infty\}$, if we recall
from~\cite{CAoRFI} that $c_{k}=\alpha_{k}-\ii\beta_{k}$, and also that
$c_{k}^{(-n)}=\alpha_{k}^{(-n)}-\ii\beta_{k}^{(-n)}$,

\noindent
\begin{eqnarray}
  \alpha_{k}-\ii\beta_{k}
  & = &
  \ii^{n}k^{n}
  \left[
    \alpha_{k}^{(-n)}-\ii\beta_{k}^{(-n)}
  \right]
  \nonumber\\
  & = &
  \ii^{n}k^{n}
  \alpha_{k}^{(-n)}
  -
  \ii^{n+1}k^{n}
  \beta_{k}^{(-n)}.
\end{eqnarray}

\noindent
We now see that the relations between the Fourier coefficients
$(\alpha_{k},\beta_{k})$ and the Fourier coefficients
$\left[\alpha_{k}^{(-n)},\beta_{k}^{(-n)}\right]$ depend on the parity of
$n$. For even $n=2j$ we have

\noindent
\begin{eqnarray}
  \alpha_{k}
  & = &
  (-1)^{j}k^{n}
  \alpha_{k}^{(-n)},
  \nonumber\\
  \beta_{k}
  & = &
  (-1)^{j}k^{n}
  \beta_{k}^{(-n)},
\end{eqnarray}

\noindent
while for odd $n=2j+1$ we have

\noindent
\begin{eqnarray}
  \alpha_{k}
  & = &
  (-1)^{j}k^{n}
  \beta_{k}^{(-n)},
  \nonumber\\
  \beta_{k}
  & = &
  (-1)^{j+1}k^{n}
  \alpha_{k}^{(-n)}.
\end{eqnarray}

\noindent
Since we have the coefficients $\alpha_{k}^{(-n)}$ and $\beta_{k}^{(-n)}$
written as integrals, we may now write $\alpha_{k}$ and $\beta_{k}$ as
integrals, first for the case of even $n=2j$,

\noindent
\begin{eqnarray}
  \alpha_{k}
  & = &
  \frac{(-1)^{j}k^{n}}{\pi}
  \int_{-\pi}^{\pi}d\theta\,
  \cos(k\theta)f_{r}^{-n\prime}(\theta),
  \nonumber\\
  \beta_{k}
  & = &
  \frac{(-1)^{j}k^{n}}{\pi}
  \int_{-\pi}^{\pi}d\theta\,
  \sin(k\theta)f_{r}^{-n\prime}(\theta),
\end{eqnarray}

\noindent
and then for the case of odd $n=2j+1$,

\noindent
\begin{eqnarray}
  \alpha_{k}
  & = &
  \frac{(-1)^{j}k^{n}}{\pi}
  \int_{-\pi}^{\pi}d\theta\,
  \sin(k\theta)f_{r}^{-n\prime}(\theta),
  \nonumber\\
  \beta_{k}
  & = &
  \frac{(-1)^{j+1}k^{n}}{\pi}
  \int_{-\pi}^{\pi}d\theta\,
  \cos(k\theta)f_{r}^{-n\prime}(\theta).
\end{eqnarray}

\noindent
Since the integrals shown are all limited as functions of $k$, given that
$f_{r}^{-n\prime}(\theta)$ is an integrable real function, once again it
is apparent that the coefficients $\alpha_{k}$ and $\beta_{k}$ which are
associated to $f_{r}(\theta)$ typically diverge as a positive power of $k$
when $k\to\infty$. Note that, when the real function $f(\theta)$ is
integrable on the whole unit circle, we are reduced to the case $n=0$, so
that the expressions for the even case reduce to the usual ones for the
Fourier coefficients.

In a previous paper~\cite{CAoRFIII} we showed that the whole Fourier
theory of integrable real functions is contained in the complex-analytic
structure introduced in~\cite{CAoRFI}. We also extended that Fourier
theory to include, not only the singular distributions discussed
in~\cite{CAoRFII}, but essentially the whole space of inner analytic
functions. We now observe that the sequences of complex Taylor
coefficients $c_{k}$ in Equation~(\ref{EQGenCoe}), which go to infinity
with $k$ not faster than a power, are exponentially bounded, according to
the definition of that term given in~\cite{CAoRFIII}, and repeated in
Equation~(\ref{EQExpLim}). This in itself suffices to show that the
corresponding power series converges to an inner analytic function, as was
shown in that paper. It also implies, as was also shown there, that the
two sequences of real coefficients $\alpha_{k}$ and $\beta_{k}$ associated
to $f(\theta)$ are both also exponentially bounded. As a consequence of
all this, the non-integrable real functions we are discussing in this
paper can be expressed as regulated Fourier series, as given in
Equation~(\ref{EQRegSer}), thus using the summation rule that was
presented in~\cite{CAoRFIII},

\begin{equation}
  f(\theta)
  =
  \frac{\alpha_{0}}{2}
  +
  \lim_{\rho\to 1_{(-)}}
  \sum_{k=1}^{\infty}
  \rho^{k}
  \left[
    \alpha_{k}\cos(k\theta)
    +
    \beta_{k}\sin(k\theta)
  \right],
\end{equation}

\noindent
which is equivalent to the fact that the non-integrable real functions
which we discussed here can be obtained as the $\rho\to 1_{(-)}$ limits of
the real parts of inner analytic functions. We see therefore that the
class of non-integrable real functions which we are examining here is also
contained in the extended Fourier theory presented in~\cite{CAoRFIII}.

\section{Conclusions and Outlook}\label{Sec06}

We have extended the close and deep relationship established in previous
papers~\cite{CAoRFI,CAoRFII}, between, on the one hand, integrable real
functions and singular Schwartz distributions, and, in the other hand,
complex analytic functions in the unit disk centered at the origin of the
complex plane, to include a large class of non-integrable real functions.
This close relationship between real functions and related objects, and
complex analytic functions, allows one to use the powerful and extremely
well-known machinery of complex analysis to deal with the real functions
and related objects in a very robust way, even if these objects are very
far from being analytic. The concept of integral-differential chains of
proper inner analytic functions, which we introduced in~\cite{CAoRFI},
played a central role in the analysis presented.

One does not usually associate non-differentiable, discontinuous and
unbounded real functions, as well as singular distributions, with single
analytic functions. Therefore, it may come as a bit of a surprise that, as
was established in~\cite{CAoRFI,CAoRFII}, essentially {\em all} integrable
real functions, as well as all {\em all} singular Schwartz distributions,
are given by the real parts of certain inner analytic functions on the
open unit disk, in the limit in which one approaches the unit circle. This
surprise is now further compounded by the fact that inner analytic
functions can represent a large class of non-integrable real functions as
well.

There are still more inner analytic functions within the open unit disk
than those that were examined here and in~\cite{CAoRFI,CAoRFII}, in
relation to integrable real functions, singular distributions and
non-integrable real functions. One important limitation in the arguments
presented here is that requiring that there be only a finite number of
non-integrable hard singularities. It may be possible, perhaps, to lift
this limitation, allowing for a denumerably infinite set of such
non-integrable singularities. It is probably not possible, however, to
allow for a densely distributed set of such singularities. Possibly, the
limitation that the number of non-integrable hard singularities be finite
may be exchanged for the limitation that the number of {\em accumulation
  points} of a denumerably infinite set of singular points with
non-integrable hard singularities be finite. We may conjecture that the
following condition might work: consider the set of all closed intervals
contained in the unit circle on which the function is integrable; consider
the point-set infinite union of all such intervals; if the resulting set
has measure $2\pi$, then it may be possible to show that the real function
can still be represented by an inner analytic function $w(z)$, and thus by
a definite set of Fourier coefficients. This condition certainly holds for
the cases discussed in this paper, and may even turn out to be the most
general possible condition leading to the results.

One interesting aspect of the work presented here is that we obtain a
definite set of Fourier coefficients associated to the non-integrable real
function $f(\theta)$. It is particularly interesting to note the fact
that, when we determine $\alpha_{0}$ during the construction of the
corresponding inner analytic function $w(z)$, we are in effect {\em
  defining}, by analytic means, the average value, over the unit circle,
of the non-integrable function $f(\theta)$, for which such a concept was
not previously defined at all. A similar remark can be made for all the
other Fourier coefficients as well. Just as was the case for integrable
real functions and singular Schwartz distributions, the non-integrable
real functions examined here can be said to be represented directly by
their sequences of Fourier coefficients.

One way to interpret the structure presented here is that, although the
non-integrable real function $f(\theta)$ is defined within separate
sections, with no predetermined relations among them, it is always
possible to define an inner analytic function that reproduces the
non-integrable real function correctly, strictly within each one of these
sections, and therefore connects them all to one another in an analytic
way, just as the interior of the unit disk connects the arcs of the unit
circle that correspond to the sections.

We believe that the results presented here enlarge the new perspective for
the analysis of real functions which was established in~\cite{CAoRFI}.
This development confirms the opinion expressed there that the use of the
theory of complex analytic functions makes it a rather deep and powerful
point of view. Since complex analysis and analytic functions constitute in
fact such a powerful tool, with so many applications in almost all areas
of mathematics and physics, it is to be hoped that further applications of
the ideas explored here will in due time present themselves.

\section*{Acknowledgments}

The author would like to thank his friend and colleague Prof. Carlos
Eugênio Imbassay Carneiro, to whom he is deeply indebted for all his
interest and help, as well as his careful reading of the manuscript and
helpful criticism regarding this work.

\bibliography{allrefs_en}\bibliographystyle{ieeetr}

\end{document}